\input ./style/arxiv-general.cfg
\documentclass[MSNbibl,number,citesort,dvips]{arxbj}
\makeatletter
   \@ifpackageloaded{graphicx}{}{\usepackage{graphicx}}
\makeatother
\usepackage{accents}%
\usepackage{mathbh}
\usepackage{mathrsfs}
\volume{22}
\issue{3}
\pubyear{2016}
\firstpage{1431}
\lastpage{1447}
\doi{10.3150/15-BEJ698}% Updated by VTEXPTS2LaTeX.exe, 29.04.2015 08:29
\docsubty{FLA}

\makeatletter
\newcommand{\eqref}[1]{(\ref{#1})}
\renewcommand{\mathring}[1]{\accentset{\circ}{#1}}
\def\overset{\mathring}
\newtheorem{theorem}{Theorem}[section]
\newtheorem{lemma}[theorem]{Lemma}
\newtheorem{corollary}[theorem]{Corollary}
\newproclaim{definition}[theorem]{Definition}
\newtheorem{proposition}[theorem]{Proposition}
\newremark{remark}[theorem]{Remark}
\newremark{example}[theorem]{Example}
\newremark{counterexample}[theorem]{Counterexample}
\def\cal{\mathcal}
\newcommand{\E}{\mathbb{E}}

\makeatother

\begin{document}
\begin{frontmatter}

\title{New results on mixture and exponential models by Orlicz spaces}

\runtitle{New results on mixture and exponential models}

\begin{aug}
\author[A]{\inits{M.}\fnms{Marina}~\snm{Santacroce}\thanksref{e1}\ead
[label=e1,mark]{marina.santacroce@polito.it}},
\author[A]{\inits{P.}\fnms{Paola}~\snm{Siri}\thanksref{e2}\ead
[label=e2,mark]{paola.siri@polito.it}}
\and
\author[A]{\inits{B.}\fnms{Barbara}~\snm{Trivellato}\thanksref{e3}\corref{}\ead[label=e3,mark]{barbara.trivellato@polito.it}}
%%%% inicialai - be tarpu
% Corresponding author: Barbara Trivellato - barbara.trivellato@polito.it% Updated by VTEXPTS2LaTeX.exe, 30.04.2015 08:30
%barbara.trivellato@polito.it% Updated by VTEXPTS2LaTeX.exe, 29.04.2015
%08:29
%\author[A]{\inits{}\fnms{}~\snm{}\corref{}\thanksref{A}
%\ead[label=e1]{}}%,
%\author[]{\inits{}\fnms{}~\snm{}\thanksref{}\ead[label=]{}}
% \and
%\author[]{\inits{}\fnms{}~\snm{}\thanksref{}\ead[label=]{}}
%%\runauthor{} %% auto
%\dedicated{}
\address[A]{Dipartimento di Scienze Matematiche ``G.L. Lagrange,''
Politecnico di Torino, Corso Duca degli Abruzzi 24, 10129 Torino, Italy.\\
\printead{e1,e2,e3}}
%\address[A]{. \printead{e1}}
%\address[]{. \printead{}}
\end{aug}

% HISTORY:
%
\received{\smonth{6} \syear{2014}}% Updated by VTEXPTS2LaTeX.exe,
%29.04.2015 08:29
%\revised{\smonth{} \syear{}}

% ABSTRACT
%
\begin{abstract}
New results and improvements in the study of nonparametric exponential
and mixture models are proposed. In particular, different equivalent
characterizations of maximal exponential models, in terms of open
exponential arcs and Orlicz spaces, are given. Our theoretical results
are supported by several examples and counterexamples and provide an
answer to some open questions in the literature.
\end{abstract}

% KEYWORDS
% visi is mazosios raides ir pagal abecele
%
\end{frontmatter}

%s1 #&#
\section{Introduction}\label{sec1}
The geometry of statistical models, called Information Geometry,
started with a paper of Radhakrishna Rao~\cite{R} and has been described in its modern
formulation by Amari \cite{A1,A2} and Amari and Nagaoka~\cite{AN}. Until the nineties, the
theory was developed only in the parametric case. The first rigorous
infinite dimensional extension has been formulated by Pistone and Sempi \cite{PS}. In
this work, the set of (strictly) positive densities has been endowed
with a structure of exponential Banach manifold, using the Orlicz space
associated to an exponentially growing Young function. The geometry of
nonparametric exponential models and its analytical properties in the
topology of the exponential Orlicz space has been also studied in
subsequent works, for example, by Gibilisco and Pistone~\cite{GP}, Pistone and Rogantin \cite{PR}, Cena and Pistone \cite{CP}.

In this paper, we develop some ideas contained in Cena and Pistone \cite{CP} and we add
several new results and improvements in the study of nonparametric
exponential and mixture models. In particular, a novelty is represented
by the introduction in this context of the time dependence, which could
allow to change the perspective from static to dynamic.

In the exponential framework, the starting point is the notion of
\textit{maximal exponential model} centered at a given positive density
$p$, introduced by Pistone and Sempi \cite{PS}. One of the main results of Cena and Pistone \cite{CP}
states that any density belonging to the maximal exponential model
centered at $p$ is connected by an \textit{open} exponential arc to $p$
and vice versa (by ``open,'' we essentially mean that the two densities
are not the extremal points of the arc). In this work, we give a proof
of this result, which is at the same time simpler and more rigorous
than the one in Cena and Pistone \cite{CP}. Moreover, we additionally prove that the
equality of the maximal exponential models centered at two (connected)
densities $p$ and $q$ is equivalent to the equality of the Orlicz
spaces referred to the same densities.
Our achievements highlight the role of the Orlicz spaces in the theory
of nonparametric exponential models and its connection with the
divergence between densities, and thus, with Information Theory. Our
theoretical results are supported by several examples and
counterexamples, which provide an answer to some open questions,
filling some gaps in the literature.

A second part of the work is devoted to the study of open mixture arcs
and contains results which are the counterpart of those obtained for
open exponential arcs. More specifically, we give the characterization
of open mixture models by establishing the equivalence between the open
mixture connection and the boundedness of the densities ratios $\frac
{q}{p}$ and $\frac{p}{q}$.

The paper is organized as follows. In Section~\ref{preliminary}, some
basic notions in the theory of Orlicz spaces are briefly recalled. The
definitions of open mixture and exponential arcs are given in Section~\ref{arcs}.
Section~\ref{mainres} contains our main results. More specifically, the
characterizations of exponential and mixture models are dealt with in
Section~\ref{char}. Densities time evolution and some geometric
properties of exponential and mixture models, namely the convexity and
the $L^1$ closure, are studied, respectively, in Sections~\ref{time}
and~\ref{geo}.

%s2 #&#
\section{Preliminaries on Orlicz spaces} \label{preliminary}

In this section, we recall some known results from the theory of Orlicz
spaces, which will be useful in the sequel.
For further details on Orlicz spaces, the reader is referred to Rao and Ren~\cite{RR1,RR2}.

Let $({\cal X}, {\mathcal F}, \mu)$ be a fixed measure space.
Young functions can be seen as generalizations of the functions
$f(x)=\frac{|x|^a}{a}$, with $a>1$, and consequently, Orlicz spaces are
generalizations of the Lebesgue spaces $L^a(\mu)$. Now, we give the
definition of Young function and of the related Orlicz space.

%de2.1 #&#
\begin{definition}
A Young function $\Phi$ is an even, convex function $\Phi\dvtx  {\mathbb
R}\rightarrow[0,+\infty]$ such that
\begin{longlist}[(iii)]
\item[(i)]
$\Phi(0)=0$,
\item[(ii)]
$\lim_{x\rightarrow\infty}\Phi(x)=+\infty$,
\item[(iii)]
$\Phi(x)<+\infty$ in a neighborhood of 0.
\end{longlist}
\end{definition}

The \textit{conjugate function} $\Psi$ of $\Phi$, is defined as $\Psi
(y)=\sup_{x\in{\mathbb R}}\{xy-\Phi(x)\}$, $\forall y \in
{\mathbb R}$ and it is itself a Young function.
From the definition of $\Psi$, the Fenchel--Young inequality
immediately follows:
%
%e1 #&#
\begin{equation}
\label{FYin} |xy|\leq\Phi(x)+\Psi(y), \qquad x,y\in{\mathbb R}.
\end{equation}

This inequality is a generalization of the classical Young inequality
$|xy|\leq\frac{|x|^a}{a}+\frac{|y|^b}{b}$ with $a,b>0$, $\frac
{1}{a}+\frac{1}{b}=1$, used in the ordinary $L^a(\mu)$ spaces.

Now, let $L^0$ denote the set of all measurable functions $u\dvtx{\cal X}
\rightarrow{\mathbb R}$ defined on $({\cal X}, {\mathcal F}, \mu)$.

%de2.2 #&#
\begin{definition}
The Orlicz space $L^\Phi(\mu)$ associated to the Young function $\Phi$
is defined as
%
%e2 #&#
\begin{equation}
L^\Phi(\mu)=\biggl\{u\in L^0 \dvt \exists \alpha>0 \mbox{ s.t. }
\int_{\cal X}\Phi(\alpha u) \,\mathrm{d}\mu<+\infty\biggr\}.
\end{equation}
\end{definition}

The Orlicz space $L^\Phi(\mu)$ is a vector space. Moreover, one can
show that it is a Banach space when endowed with the \textit{Luxembourg norm}
%
%e3 #&#
\begin{equation}
\|u\|_{\Phi,\mu}=\inf\biggl\{k>0 \dvt \int_{\cal X}\Phi
\biggl(\frac{u} k \biggr) \,\mathrm{d}\mu\leq1\biggr\}.
\end{equation}

Consider the Orlicz space $L^\Phi(\mu)$ with the Luxembourg norm $\|
\cdot\|_{\Phi,\mu}$ and denote by $B(0,1)$ the open unit ball and by
$\overline{B(0,1)}$ the closed one. Let us observe that
\begin{eqnarray*}
u&\in& B(0,1) \quad\iff \quad\exists \alpha>1 \mbox{ s.t. } \int_{\cal X}
\Phi (\alpha u) \,\mathrm{d}\mu\leq1,
\\
u&\in&\overline{B(0,1)} \quad\iff\quad \int_{\cal X} \Phi(u) \,\mathrm{d}\mu\leq1.
\end{eqnarray*}

Moreover, the Luxembourg norm is equivalent to the \textit{Orlicz norm}
%
%e4 #&#
\begin{equation}
N_{\Phi,\mu}(u)=\sup_{v\in L^\Psi(\mu)\dvt \int_{\cal X}\Psi(v) \,\mathrm{d}\mu
\leq1}\biggl\{\int
_{\cal X}|uv| \,\mathrm{d}\mu \biggr\},
\end{equation}
where $\Psi$ is the conjugate function of $\Phi$.

It is worth to recall that the same Orlicz space can be related to
different \textit{equivalent} Young functions.

%de2.3 #&#
\begin{definition}
Two Young functions $\Phi$ and $\Phi'$ are said to be equivalent if
there exists $x_0>0$, and two positive constants $c_1<c_2$ such that,
$\forall x\geq x_0$,
\[
\Phi(c_1 x)\leq\Phi'(x)\leq\Phi(c_2 x).
\]
\end{definition}

 In such a case the Orlicz spaces $L^{\Phi} (\mu)$ and $L^{\Phi
'} (\mu)$ are equal as sets and have equivalent norms as Banach spaces.

From now on, we consider a probability space $({\cal X}, {\mathcal F},
\mu)$ and we denote with $\cal P$ the set of all densities which are
positive $\mu$-a.s.
Moreover, we use $\E_p$ to denote the integral with respect to $p\,\mathrm{d}\mu$,
for each fixed $p\in\cal P$.

In the sequel, we use the Young function $\Phi_1(x)= \cosh(x)-1$, which
is equivalent to the more commonly used $\Phi_2(x)= \mathrm{e}^{|x|}-|x|-1$.

We recall that the conjugate function of $\Phi_1(x)$ is $\Psi_1(y)= \int_0^y\sinh^{-1}(t) \,\mathrm{d}t$, which, in its turn, is equivalent to $\Psi_2(y)=
(1+|y|)\log(1+|y|)-|y|$.

Finally, in order to stress that we are working with densities $p\in
\cal P$, we will denote with $L^{\Phi_1} (p)$ the Orlicz space
associated to $\Phi_1$, defined with respect to the measure induced by
$p$, that is,
%
%e5 #&#
\begin{equation}
L^{\Phi_1} (p)=\bigl\{u \in L^0 \dvt \exists \alpha>0 \mbox{ s.t. } \E
_p\bigl(\Phi_1(\alpha u)\bigr)<+\infty\bigr\}.
\end{equation}

It is worth to note that, in order to prove that a random variable $u$
belongs to $L^{\Phi_1} (p)$, it is sufficient to check that
$\E_p(\mathrm{e}^{\alpha u})<+\infty$, with $\alpha$ belonging to an open
interval containing 0.

%s3 #&#
\section{Mixture and exponential arcs}\label{arcs}

In this section, we recall the definitions of mixture and exponential
arcs, and some related results.

%de3.1 #&#
\begin{definition}\label{mix}
Two densities $p, q \in{\cal P}$ are connected by an open mixture arc
if there exists an open interval $I \supset[0,1]$ such that $p(\theta
)=(1-\theta)p+\theta q$ belongs to $\cal P$, for every $\theta\in I$.
\end{definition}

%de3.2 #&#
\begin{definition}\label{exp}
Two densities $p, q \in{\cal P}$ are connected by an open exponential
arc if there exists an open interval $I \supset[0,1]$ such that
$p(\theta)\propto p^{(1-\theta)}q^{\theta}$ belongs to $\cal P$, for
every $\theta\in I$.
\end{definition}

In the following proposition, we give an equivalent definition of
exponential connection by arcs.

%pr3.3 #&#
\begin{proposition}
$p, q \in{\cal P}$ are connected by an open exponential arc iff there
exist an open interval $I \supset[0,1]$ and a random variable $u\in
L^{\Phi_1}(p)$, such that
$p(\theta)\propto \mathrm{e}^{\theta u}p$ belongs to $\cal P$, for every $\theta
\in I$ and $p(0)=p,  p(1)=q$.
\end{proposition}

\begin{pf}
Let us assume that $p, q \in{\cal P}$ are connected by an open
exponential arc, that is, $\int_{\cal X} p^{(1-\theta)}q^{\theta} \,\mathrm{d}\mu
<+\infty$, for any $\theta\in I$.
Since
\[
\int_{\cal X} p^{(1-\theta)}q^{\theta} \,\mathrm{d}\mu=
\E_p\biggl(\biggl(\frac{q}p \biggr)^\theta\biggr)=
\E_p \bigl(\mathrm{e}^{\theta u}\bigr) \qquad\mbox{with } u=\log
\frac{q} p,
\]
then $u\in L^{\Phi_1}(p)$. Moreover $p(\theta)\propto \mathrm{e}^{\theta u}p$
belongs to $\cal P$, for every $\theta\in I$ and $p(0)=p,  p(1)=q$.

The converse follows immediately, observing that $q=p(1)\propto\mathrm{e}^up$,
that
is, $u=\log\frac{q} p+c$.
\end{pf}

The connections by open mixture arcs and by open exponential arcs are
equivalence relations (see Cena and Pistone \cite{CP} for the proofs).

In the following, we recall the definition of the cumulant generating
functional and its properties, in order to introduce the notion of
maximal exponential model.
In the next section, the maximal exponential model at $p$ is proved to
coincide with the set of all densities $q \in{\cal P}$ which are
connected to $p$ by an open exponential arc.

Let us denote
\[
L^{\Phi_1}_{0}(p)=\bigl\{u \in L^{\Phi_1}(p)\dvt
\E_p(u)=0\bigr\}.
\]

%de3.4 #&#
\begin{definition}\label{K}
The cumulant generating functional is the map
%
%e6 #&#
\begin{eqnarray}
K_p\dvtx && L^{\Phi_1}_{0}(p) \longrightarrow
[0,+\infty],
\nonumber
\\[-8pt]
\\[-8pt]
\nonumber
&& u \longmapsto  \log\E_p \bigl(\mathrm{e}^u\bigr).
\end{eqnarray}
\end{definition}

%th3.5 #&#
\begin{theorem}\label{ThPS}
The cumulant generating functional $K_p$ satisfies the following properties:
\begin{longlist}[(ii)]
\item[(i)] $K_p(0)=0$; for each $u\neq0$, $K_p(u)>0$.
\item[(ii)] $K_p$ is convex and lower semicontinuous, moreover its
proper domain
\[
\operatorname{dom} K_p=\bigl\{u \in L^{\Phi_1}_{0}(p)
\dvt K_p(u)< +\infty\bigr\}
\]
is a convex set which contains the open unit ball of $L^{\Phi
_1}_{0}(p)$. In particular, its interior $\overset{\operatorname{dom}
K_p}$ is a nonempty convex set.
\end{longlist}
\end{theorem}

For the proof, one can see Pistone and Sempi \cite{PS}.

%de3.6 #&#
\begin{definition}
For every density $p \in{\cal P}$, the maximal exponential model at
$p$ is
\[
{\cal E}(p)=\bigl\{q=\mathrm{e}^{u-K_p(u)} p \dvt u\in\overset {
\operatorname{dom} K_p} \bigr\} \subseteq{\cal P}.
\]
\end{definition}

%re3.7 #&#
\begin{remark}
We have defined $K_p$ on the set $L^{\Phi_1}_{0}(p)$ because centered
random variables guarantee the uniqueness of the representation of
$q\in{\cal E}(p)$.
\end{remark}

%s4 #&#
\section{Main results on mixture and exponential models}\label{mainres}

%s4.1 #&#
\subsection{Characterizations}\label{char}

In the sequel, we use the notation $D(q\|p)$ to indicate the
Kullback--Leibler divergence of $q\cdot\mu$ with respect to $p\cdot\mu
$ and we simply refer to it as the divergence of $q$ from $p$.

We first state two results related to Orlicz spaces, which will be used
in the sequel.
Their proofs can be found in Cena and Pistone \cite{CP}.

%pr4.1 #&#
\begin{proposition}\label{L1}
Let $p$ and $q$ belong to $\mathcal P$ and let $\Phi$ be a Young
function.

The Orlicz spaces $L^\Phi(p)$ and $L^\Phi(q)$ coincide if
and only if their norms are equivalent.
\end{proposition}

%le4.2 #&#
\begin{lemma}\label{finitediv}
Let $p,q \in\mathcal P$, then
$D(q\|p)<+\infty \Longleftrightarrow \frac{q}{p}\in L^{\Psi
_1}(p) \Longleftrightarrow \log\frac{q}{p}\in L^{1}(q)$.
\end{lemma}

From Lemma~\ref{finitediv}, we can prove the following result.

%th4.3 #&#
\begin{theorem}\label{immersion}
Let $p,q \in\mathcal P$.
If $D(q\|p)<+\infty$ then $L^{\Phi_1}(p)\subseteq L^{1}(q)$.
\end{theorem}

\begin{pf}
Let us consider $u\in L^{\Phi_1}(p)$, i.e. $\E_p(\Phi_1(\alpha
u))<+\infty$ for some $\alpha>0$.

Note that, by Lemma~\ref{finitediv}, the hypothesis $D(q\|p)<+\infty$
is equivalent to $\frac{q}{p}\in L^{\Psi_1}(p)$, i.e.
$\E_p(\Psi_1(\beta\frac{q} p))<+\infty$ for some
$\beta>0$.

Thus, using the Fenchel--Young inequality \eqref{FYin}
and taking the expectation, we deduce that
\[
\alpha\beta\E_q\bigl(|u|\bigr)= \alpha\beta\E_p\biggl(|u|
\frac{q} p\biggr)\leq\E _p\bigl(\Phi_1(\alpha u)
\bigr)+\E_p\biggl(\Psi_1\biggl(\beta\frac{q}p
\biggr)\biggr)<+\infty.
\]
\upqed\end{pf}

%re4.4 #&#
\begin{remark}
In Theorem~\ref{immersion}, $\E_p (\Psi_1 (\beta\frac{q} p))<+\infty$ for some $\beta>0$ equals
$\E_p (\Psi_1 ( \frac{q} p))$ $<+\infty$. In
fact, ${\Psi_1}$ is equivalent to $\Psi_2$ and it is easy to check that
$\Psi_2$ satisfies the generalized $\Delta_2$ condition
\[
\Psi_2(\beta y)\le\max{\bigl(\beta^2,1\bigr)}
\Psi_2(y).
\]
Assume $y>0$ and observe that $\Psi_2(y)= (1+|y|)\log(1+|y|)-|y|$
admits the representation
\[
\psi_2(y)=\int_0^y
\frac{y-\tau}{1+\tau}\,\mathrm{d}\tau.
\]
Therefore,
\[
\psi_2(\beta y)=\beta^2 \int_0^y
\frac{y-\tau}{1+\beta\tau}\,\mathrm{d}\tau\le \max{\bigl(\beta^2,1\bigr)}
\Psi_2(y).
\]
\end{remark}

%%%%%%%%%%%%%
%Figura 1
%f1 #&#
\begin{figure}[b]

\includegraphics{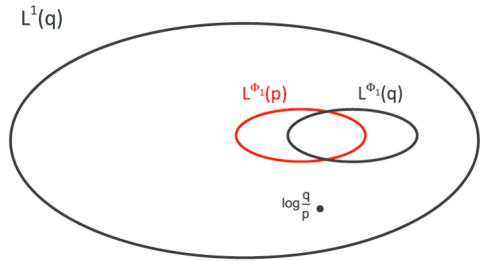}

\caption{The case when $\log\frac{q}{p}\in L^{1}(q)$, i.e. $D(q\|
p)<+\infty$.}\vspace*{5pt}
\label{Fig1}
\end{figure}

%%%%%%%%%%%%%

%%%%%%%%%%%%%
%Figura 2
%f2 #&#
\begin{figure}

\includegraphics{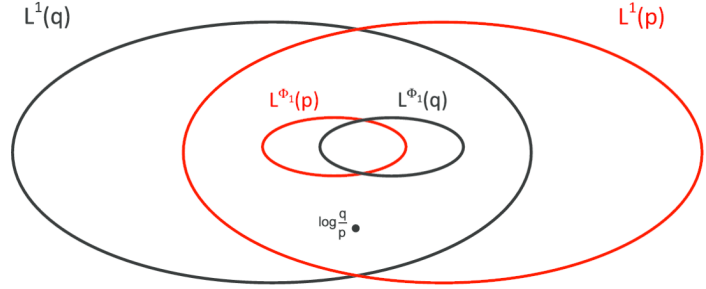}

\caption{The case when $\log\frac{q}{p}\in L^{1}(q)\cap L^{1}(p)$, i.e.
$D(q\|p)<+\infty$ and $D(p\|q)<+\infty$.}
\label{Fig2}
\end{figure}

Figures~\ref{Fig1} and \ref{Fig2} show the geometry described in
Theorem~\ref{immersion}.\

The next two results are technical preliminaries to Theorem~\ref{main}.
In particular, Proposition~\ref{P1} gives a sufficient condition on
Orlicz norms in order to have $\log\frac{q}{p}\in L^{1}(p)$.

%le4.5 #&#
\begin{lemma}\label{L2}
Let $p$ and $q$ belong to $\mathcal P$ and let $M$ be any positive
constant. Then
\[
\biggl\|\mathbh{1}_{({q}/{p}>M)}\log\frac{q}{p}\biggr\|_{\Phi_1,p}<+\infty.
\]
\end{lemma}
\begin{pf}
Let us denote $A = \{\frac{q}{p}>M\}$ and recall that $\mathbh{1}_A \log
\frac{q}{p}  \in  L^{\Phi_1}(p)$ if and only if $\E_p(\mathrm{e}^{\alpha\mathbh{1}_A\log{q}/{p}})$ $<+\infty$ for any $\alpha\in(-\varepsilon,
+\varepsilon)$ with $\varepsilon$ sufficiently small. Since
\[
\label{eq1}\E_p\bigl(\mathrm{e}^{\alpha\mathbh{1}_A\log({q}/{p})}\bigr)\le 1+
\E_p\biggl(\mathbh{1}_A\biggl(\frac{q}{p}
\biggr)^\alpha\biggr),
\]
when $0<\alpha<\varepsilon<1$, by Jensen inequality
\[
\E_p\biggl(\mathbh{1}_A\biggl(\frac{q}{p}
\biggr)^\alpha\biggr)\le\E_p \biggl(\biggl(\frac{q}{p}
\biggr)^\alpha\biggr)\le1,
\]
while, when $-\varepsilon<\alpha<0$,
\[
\E_p\biggl(\mathbh{1}_A\biggl(\frac{q}{p}
\biggr)^\alpha\biggr)\le M^\alpha.
\]
\upqed\end{pf}

%pr4.6 #&#
\begin{proposition}\label{P1}
If $\|\cdot\|_{\Phi_1,p}\le c \|\cdot\|_{\Phi_1,q}$, then
$\log\frac{q}{p} \in L^{\Phi_1}(p)$.
\end{proposition}

\begin{pf}
First, we write
\begin{eqnarray*}
\log\frac{q}{p}&=&\mathbh{1}_{({q}/{p}>M)}\log\frac{q}{p}+
\mathbh{1}_{({q}/{p}\leq M)}\log\frac{q}{p}
\\
&=&\mathbh{1}_{({q}/{p}>M)}\log\frac{q}{p}-\mathbh{1}_{({p}/{q}\geq M^{-1})}\log
\frac{p}{q}.
\end{eqnarray*}
By hypothesis, we have $\|\cdot\|_{\Phi_1,p}\le c \|\cdot\|_{\Phi
_1,q}$. Therefore,
\[
\biggl\|\log\frac{q}{p}\biggr\|_{\Phi_1,p}\le\biggl\|\mathbh{1}_{({q}/{p}>M)}\log
\frac
{q}{p}\biggr\|_{\Phi_1,p}+ c\biggl\|\mathbh{1}_{({p}/{q}\geq M^{-1})}
\log\frac{p}{q}\biggr\|
_{\Phi_1,q}.
\]
%
%Now, the conclusion immediately follows from Lemma~\ref{L2}.
Now, the conclusion immediately follows from Lemma \ref{L2},
noting that its result
holds also when the strict inequality $\frac{q}{p} > M$ is replaced by
$\frac{q}{p}\geq M$.
\end{pf}

The following theorem is an important improvement of Theorem~21 of Cena and Pistone \cite
{CP}. In particular, the main point is the equivalence between the
equality of the exponential models $\mathcal E(p)$ and $\mathcal E(q)$
and the equality of the Orlicz spaces $L^{\Phi_1}(p)$ and $L^{\Phi
_1}(q)$. Moreover, we show that if a density belongs to the maximal
exponential model at $p$, there exists an open exponential arc
connecting the two densities and vice versa.
This result was first stated in Theorem~21 of Cena and Pistone \cite{CP}. However, the
proof is somehow involved and imprecise in some steps. Here, we give a
simpler and rigorous one.\vadjust{\goodbreak}

%th4.7 #&#
\begin{theorem}\label{main}
Let $p,q \in\mathcal P$. The following statements are equivalent:
\begin{longlist}[(iii)]
\item[(i)] $q\in\mathcal E(p)$;
\item[(ii)] $q$ is connected to $p$ by an open exponential arc;
\item[(iii)] $\mathcal E(p)=\mathcal E(q)$;
\item[(iv)] $L^{\Phi_1}(p)=L^{\Phi_1}(q)$;
\item[(v)] $\log\frac{q}{p} \in L^{\Phi_1}(p)\cap L^{\Phi_1}(q)$;
\item[(vi)] $\frac{q}{p}\in L^{1+\varepsilon}(p)$  and  $\frac{p}{q}
\in L^{1+\varepsilon}(q)$, for some $\varepsilon>0$.
\end{longlist}
\end{theorem}
\begin{pf}
We first show the equivalence of the first two statements.
If $q \in\mathcal E(p)$, $q \propto \mathrm{e}^u p$ for some $u\in\overset
{\operatorname{dom}  K_p}$. Since also $0\in\overset{\operatorname
{dom}  K_p}$ and
$\overset{\operatorname{dom}  K_p}$ is an open convex set, we deduce
that $p(\theta)\propto \mathrm{e}^{\theta u}p$
is an open exponential arc containing $p$ and $q$, for $\theta$ in an
open interval $I\supset[0,1]$.

Vice versa, assume $q$ is connected to
$p$ by an open exponential arc $p(\theta)\propto \mathrm{e}^{\theta u}p$ with
$q=p(1)$. Since the exponential arc is defined for $\theta$ in an open
interval $I\supset[0,1]$, we can always choose $\overline{\theta} u\in
\operatorname{dom}  K_p$, with $\overline{\theta}>1$. We conclude that $q\in
\mathcal E(p)$, by observing that $u$ is a convex combination of
$\overline{\theta} u$ and $0\in\overset{\operatorname{dom}  K_p,}$ and
thus, it belongs to $\overset{\operatorname{dom}  K_p}$.

Note that (iii) immediately implies (i). On the other hand, if we
assume there exists an open exponential arc connecting $p$ and $q$, the
equality of the exponential models (iii) follows from the fact that
the connection through open exponential arcs is an equivalence
relation.

%Finally, let us show the equivalence of the previous statements with
%(iv) and (v).
The equivalence of the previous statements with (v) is clearly
proved in Theorem 21 of \cite{CP}. Let us show the equivalence with (iv).
(ii)${}\Rightarrow{}$(iv) is proved in Theorem~19 of Cena and Pistone~\cite{CP}.

Conversely, if we assume $L^{\Phi_1}(p)=L^{\Phi_1}(q)$, by
Propositions \ref{L1} and \ref{P1}, $\log\frac{q}{p} \in
L^{\Phi_1}(p)=L^{\Phi_1}(q)$. Denoting $u=\log\frac{q}{p}$, one can
observe that, for any $\theta\in(-\varepsilon,+\varepsilon)$ with $\varepsilon
>0$, $\E_p(\mathrm{e}^{\theta u})<+\infty$ and $\E_q(\mathrm{e}^{\theta u})=\E
_p(\mathrm{e}^{(1+\theta) u})<+\infty$.
Therefore, using also Jensen inequality, one concludes that there
exists an open exponential arc connecting $p$ and $q$.

%%%%%%%%%%%%% %%%%%%%%%%%%%
%Figura 3
%f3 #&#
\begin{figure}

\includegraphics{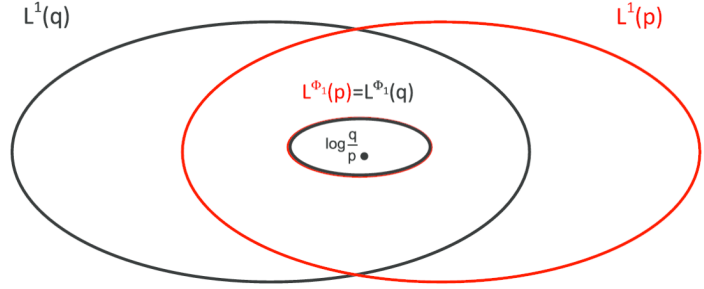}

\caption{The particular case of Figure \protect\ref{Fig2} when $\log\frac{q}{p}
\in L^{\Phi_1}(p)\cap L^{\Phi_1}(q)$ i.e. $q \in\mathcal{E}(p)$.}
\label{Fig3}
\end{figure}

The equivalence between (v) and (vi) easily follows from the
definition of Orlicz spaces.
\end{pf}

The geometry of Theorem \ref{main} is shown in Figure \ref{Fig3}.\vadjust{\goodbreak}

%co4.8 #&#
\begin{corollary}
If $q \in\mathcal{E}(p)$, then the divergences $D(q\|p)<+\infty$ and
$D(p\|q)<+\infty$.
\end{corollary}
\begin{pf}
The thesis follows from (vi) of Theorem~\ref{main}, taking into
account that $\frac{q}{p}\in L^{1+\varepsilon}(p)$ is a sufficient
condition for $D(q\|p)<+\infty$.
\end{pf}

%%%%%%%%%%%%%

The converse of this corollary does not hold. In the following, we
provide a counterexample, which at the same time, answers to an open
question raised by Cena and Pistone \cite{CP}.

%co4.9 #&#
\begin{counterexample}\label{divergenza}
Let us denote by ${\cal X}=[0,1]$, ${\mathcal F}={\cal B}([0,1])$ and
$\mu$ the corresponding Lebesgue measure.
We consider the trivial density $p(x)=1$ and
%
%e7 #&#
\begin{equation}
q(x)= C \sum_{n=1}^\infty\frac{1} {n^3 C_n}
\biggl(x-\biggl(1 -\frac{1} {n}\biggr)\biggr)^{-{n} /{(n+1)}}
 \mathbh{1}_{( 1 -{1}/ {n}, 1 -{1}/ {(n+1)}]}(x),
\end{equation}
where
%
%e8 #&#
%e9 #&#
\begin{eqnarray}
C_n&=&\int_{1 -{1}/ {n}}^{1 -{1}/ {(n+1)}}\biggl(x-\biggl(1
-\frac{1}{n}\biggr)\biggr)^{-{n}/ {(n+1)}} \,\mathrm{d}x= \frac{n+1}{\sqrt[n+1]{n(n+1)}},
\\
C&=&\Biggl( \sum_{n=1}^\infty
\frac{1} {n^3}\Biggr)^{-1}.
\end{eqnarray}

Now we prove that, $\forall\varepsilon>0$,
\begin{eqnarray*}
\E_{\mu}\bigl(q^{1+\varepsilon}\bigr) &=& \int_0^1
q(x)^{1+\varepsilon} \,\mathrm{d}x = \sum_{n=1}^\infty
\biggl(\frac{C} {n^3 C_n} \biggr)^{1+\varepsilon}
\int_{1 -{1}/ {n}}^{1-{1}/ {(n+1)}}
{\biggl(x- \biggl(1 -\frac{1} {n}\biggr)\biggr)
^{-{n(1+\varepsilon)} /{(n+1)}}} \,\mathrm{d}x\\
&=&+\infty.
\end{eqnarray*}
In fact,
\[
\int_{1 -{1}/ {n}}^{ 1-{1} /{(n+1)}} {\biggl(x-\biggl(1-
\frac{1}{n}\biggr)\biggr)^{-{n(1+\varepsilon)}/ {(n+1)}}} \,\mathrm{d}x
\]
converges when $\frac{n(1+\varepsilon)} {n+1}<1$, that is when $n \varepsilon<1$.
Then, for any choice of $ \varepsilon>0$, we can find infinitely many
$n>\frac{1} \varepsilon$, such that the integral above does not converge
and, as a consequence, $\E_{\mu}(q^{1+\varepsilon})=+\infty$.
Therefore, by (vi) of Theorem~\ref{main}, we deduce that $q\notin{\cal E}(1)$.

On the other hand, it can be proved that both $D(q\|p)<+\infty$ and
$D(p\|q)<+\infty$ and this concludes the counterexample.

We prove the last statements in the \hyperref[app]{Appendix}.
\end{counterexample}

%re4.10 #&#
\begin{remark}
From a geometric point of view, equality $L^{\Phi_1}(p)=L^{\Phi
_1}(q)$ in Theorem~\ref{main} is important. On the one hand, it implies
that the \textit{exponential transport mapping}, or \textit
{e-transport}, $ ^e\mathbb{U}^q_p\dvtx u \rightarrow u- \mathbb E_q(u)$ from
$L_0^{\Phi_1}(p)$ to $L_0^{\Phi_1}(q)$ is well defined. On the other
hand, due to Proposition~22 of Cena and Pistone \cite{CP}, it also implies that $L^{\Psi
_1}(p)= \frac{q}p L^{\Psi_1}(q)$.\vspace*{-1pt} As a consequence, the \textit{mixture
transport mapping}, or \textit{m-transport}, $ ^m\mathbb{U}^q_p\dvtx v
\rightarrow\frac{p}q v$ from $L_0^{\Psi_1}(p)$ to $L_0^{\Psi_1}(q)$ is
well-defined. For details on the applications of nonparametric
information geometry to statistical physics, see Pistone \cite{P}.
\end{remark}

The next theorem is the counterpart of Theorem~\ref{main} for open
mixture arcs. One of the equivalences is an improvement of Proposition~15 in Cena and Pistone \cite{CP}. In Cena and Pistone \cite{CP}, it is shown that densities connected by
open mixture arcs have bounded away from zero ratios. Here, we
additionally show the converse implication, giving a characterization
of open mixture models. Moreover, one can see that the key role for
being connected by open mixture either exponential arcs is played by
the ratios $\frac{q}{p}$ and $\frac{p}{q}$ which have to be bounded or
integrable in some sense.

Given $p\in{\cal P}$, we denote by ${\cal M}(p)$ the set of all
densities $q \in{\cal P}$ which are connected to $p$ by an open
mixture arc.

%th4.11 #&#
\begin{theorem}\label{main2}
Let $p,q \in\mathcal P$. The following statements are equivalent:
\begin{longlist}[(iii)]
\item[(i)] $q\in{\cal M}(p)$;
\item[(ii)] $\mathcal M(p)=\mathcal M(q)$;
\item[(iii)] $\frac{q}{p}, \frac{p}{q} \in L^{\infty}$.
\end{longlist}
\end{theorem}
\begin{pf}
The equivalence between (i) and (ii) follows since the relation of
connection through open mixture arcs is an equivalence relation.

Now we show that $p, q\in\mathcal P$ are connected by open mixture
arcs if and only if
\[
c_1< \frac{q}{p}< c_2 \qquad\mbox{with }
0<c_1<1< c_2.
\]
Assume $p$ and $q$ are connected by an open mixture arc that is
$p(\lambda)=\lambda q+(1-\lambda)p$ belongs to $\mathcal P$ for all
$\lambda\in(-\alpha, 1+\beta)\supset[-\varepsilon, 1+\varepsilon]$ with
$\varepsilon>0$.
Since $p(-\varepsilon)  \mbox{ and }  p(1+\varepsilon)\in\mathcal P$, it is
easy to see that $\frac{\varepsilon}{1+\varepsilon} <\frac{q}{p}<\frac
{1+\varepsilon}{\varepsilon}$.

To check the other implication, one observes that $p(\lambda)=\lambda
q+(1-\lambda)p$ belongs to $\mathcal P$ for any $\lambda\in(\frac{1}{1-c_2},\frac{1}{1-c_1})$.
\end{pf}

%pr4.12 #&#
\begin{proposition}\label{cor}
Let $p,q \in\mathcal P$. If $p$ and $q$ are connected by an open
mixture arc, then they are also connected by an open exponential arc.
\end{proposition}
\begin{pf}
The result immediately follows from Theorems \ref{main} and \ref{main2}.
\end{pf}

The converse implication does not hold, as the following counterexample shows.

%co4.13 #&#
\begin{counterexample}\label{beta1}
Consider the family of beta densities $p(\beta)\propto x^{\beta-1}$,
$x\in[0,1]$, with $\beta\in(0, +\infty)$.

It is easy to see that given two densities $p=p(\beta_1)$ and $q=p(\beta
_2)$, $\beta_1<\beta_2$, they are connected by an open exponential arc
but they are not connected by an open mixture arc. In fact, on the one
hand the open exponential arc
$p(\theta)\propto p^{1-\theta}q^{\theta}\propto x^{(1-\theta)\beta
_1+\theta\beta_2-1}$ is still in the family of beta densities for any
$\theta\in(-\frac{\beta_1}{\beta_2-\beta_1},+\infty)\supset[0,1]$.\vspace*{1pt}

On the other hand, by Theorem~\ref{main2}, it does not exist an open
mixture arc connecting $p$ and~$q$, since the ratio $\frac{q}{p}\propto
x^{\beta_2-\beta_1}$ is not bounded below by a positive constant.
\end{counterexample}

%re4.14 #&#
\begin{remark}
In Proposition~15 of Cena and Pistone \cite{CP}, by different arguments, it is
shown that if $p$ and $q$ are connected by an open mixture arc, then
the Orlicz spaces $L^{\Phi}(p)$ and $L^{\Phi}(q)$ coincide for any
Young function $\Phi$. When $\Phi=\Phi_1$ the result immediately
follows from Theorems~\ref{main} and \ref{main2}.
\end{remark}

%s4.2 #&#
\subsection{Exponential models and densities time evolution}\label{time}

In this paragraph, we introduce the time perspective in the study of
exponential models. As far as we are aware, this is the first attempt
in this direction.

Let us consider a filtration $\mathscr F=\{\mathcal F_t\dvt t\in[0,T]\}$
on the probability space $({\cal X}, {\mathcal F}, \mu)$ such that
$\mathcal F=\mathcal F_T$. Let $p\in{\cal P}$ and denote by $p_t=\E_\mu
(p|\mathcal F_t)$.

The following proposition gives a condition for the exponential
connection stability of the restrictions $p_t$ over time. From a
geometrical point of view, this result means that divergence finiteness
is preserved.

%pr4.15 #&#
\begin{proposition}\label{Ept1}
Let $t_1, t_2 \in[0,T]$, $t_1\leq t_2$.
If $p_{t_2}\in{\cal E}(p_{t_1})$ then $p_{s}\in{\cal E}(p_{t_1})$,
$\forall t_1\leq s < t_2$.
\end{proposition}

\begin{pf}
Let $ t_1\leq s < t_2$.
From condition (vi) of Theorem~\ref{main}, it is enough to prove that $\E
_{p_{t_1}}((\frac{p_s}{p_{t_1}})^{1+\varepsilon}
)<+\infty$ and $\E_{p_{s}}((\frac{p_{t_1}}{p_s}
)^{1+\varepsilon})<+\infty$, starting from the hypothesis that $\E
_{p_{t_1}}((\frac{p_{t_2}}{p_{t_1}})^{1+\varepsilon}
)<+\infty$ and $\E_{p_{t_2}}((\frac{p_{t_1}}{p_{t_2}}
)^{1+\varepsilon})<+\infty$.

 Since $x^{1+\varepsilon}$ is a convex function, by Jensen
inequality we get
%
%e10 #&#
\begin{equation}
\label{dyn} \biggl(\frac{p_s}{p_{t_1}}\biggr)^{1+\varepsilon}=\bigl(
\E_\mu (p_{t_2}|\mathcal F_s)
\bigr)^{1+\varepsilon}\frac{1}{p_{t_1}^{1+\varepsilon
}}\leq\E_\mu\bigl(p_{t_2}^{1+\varepsilon}|
\mathcal F_s\bigr)\frac{1}{p_{t_1}^{1+\varepsilon}}.
\end{equation}
Since $t_1\le s$ we have $\E_\mu(\cdot|\mathcal F_s)=\E
_{p_{t_1}}(\cdot|\mathcal F_s)$. Then, taking the
expectation with respect to $p_{t_1}$ in~\eqref{dyn}, we deduce that
\[
\E_{p_{t_1}}\biggl(\biggl(\frac{p_s}{p_{t_1}}\biggr)^{1+\varepsilon} \biggr)
\leq\E_{p_{t_1}}\biggl(\E_{p_{t_1}}\bigl(p_{t_2}^{1+\varepsilon}|
\mathcal F_s\bigr)\frac{1}{p_{t_1}^{1+\varepsilon}}\biggr)\le\E_{p_{t_1}}
\biggl( \biggl(\frac{p_{t_2}}{p_{t_1}}\biggr)^{1+\varepsilon}\biggr)<+\infty.
\]

Since also $x^{-(1+\varepsilon)}$ is a convex function, the other
condition follows in a similar way.
\end{pf}

%re4.16 #&#
\begin{remark}
As a consequence of the previous proposition, it is straightforward to
observe that if $p_{s_0}\notin{\cal E}(p_{t_1})$ for some $s_0\geq
t_1$, then $p_{s}\notin{\cal E}(p_{t_1})$
$\forall s_0\leq s \leq T$.
\end{remark}

From the above proposition, we immediately get the following result.

%co4.17 #&#
\begin{corollary}\label{Exp1}
If $p\in{\cal E}(1)$ then $p_{s}\in{\cal E}(1)$,
$\forall 0\leq s < T$.
\end{corollary}

As an application of this corollary, we can use the family of beta
densities introduced in the previous paragraph to give a concrete
example of a density belonging to ${\cal E}(1)$, along with its restrictions.

%ex4.18 #&#
\begin{example}\label{beta}
Let us denote by ${\cal X}=[0,1]$, ${\mathcal F}={\cal B}([0,1])$ and
$\mu$ the corresponding Lebesgue measure.
Define the filtration $\mathscr F=\{\mathcal F_t\dvt t\in[0,T]\}$ on
${\cal X}$, by choosing ${\mathcal F}_t = \sigma( [0,s] \dvt 0\leq s \leq
t )$.
Let $p$ be any density on $({\cal X}, {\mathcal F}, \mu)$.
Then, due to the particular choice of the filtration, the restriction
$p_t=\E_\mu(p|\mathcal F_t)$ can be written as\vspace*{-1pt}
%
%e11 #&#
\begin{equation}
\label{pt} p_t(x)=p(x) \mathbh{1}_{[0,t]}(x)+\frac{1-F(t)}{1-t}
\mathbh{1}_{(t,1]}(x),
\end{equation}
where $F(t)=\int_0^t p(x)   \,\mathrm{d}x, \forall t\in[0,1]$.

It is worth noting that $p_1=p$ and $p_0=1$ a.s.

Let us now fix $p(x)=\beta x^{\beta-1}$, with $\beta>0$.
It is easy to find an $\varepsilon>0$ such that\vspace*{-1pt}
\begin{eqnarray*}
\E_{\mu}\bigl(p^{1+\varepsilon}\bigr)&=&\int_0^1
\beta^{1+\varepsilon} x^{(\beta-1)(1+\varepsilon)} \,\mathrm{d}x<+\infty,
\\[-1pt]
\E_{\mu}\bigl(p^{-\varepsilon}\bigr)&=&\int_0^1
\beta^{-\varepsilon} x^{(\beta
-1)(-\varepsilon)} \,\mathrm{d}x<+\infty.
\end{eqnarray*}
So, from (vi) of Theorem~\ref{main}, we can conclude that $p\in{\cal E}(1)$.

 With this choice of $p$, by \eqref{pt},\vspace*{-1pt}
\[
p_t(x)=\beta x^{\beta-1} \mathbh{1}_{[0,t]}(x)+
\frac{1-t^\beta}{1-t} \mathbh{1}_{(t,1]}(x) \qquad\forall t \in[0,1],\vspace*{-1pt}
\]
and we can prove that $p_t\in{\cal E}(1)$, in a similar way.
\end{example}

In general, the converse of Proposition~\ref{Ept1} does not hold and
below we give a counterexample with $t_1=0$ and $t_2=\frac{1} 2$.

Using the same counterexample, we also define a density $q\in{\cal
E}(1)$ (along with its restrictions) such that, $\forall t\leq t_0$,
$p_t=\E_\mu(q|\mathcal F_t)$, for a fixed $t_0<\frac{1}2 $.

%co4.19 #&#
\begin{counterexample}\label{co4.19}
Using the same filtered probability space as in Example~\ref{beta},
we consider the density
%
%e12 #&#
\begin{eqnarray}
p(x)&=& C\sum_{n=1}^\infty\frac{1} {n^3 C_n}
\biggl[ \biggl(x-\biggl(\frac{1} 2 -\frac{1} {2n}\biggr)
\biggr)^{-{n}/{(n+1)}} \mathbh{1}_{({1}/ 2 -{1}/ {(2n)}, {1}/ 2 -{1}/{(2(n+1))}]}(x)
\nonumber
\\[-8pt]
\\[-8pt]
\nonumber
& &\hspace*{54pt}{}+ \biggl(\biggl(\frac{1} 2 +\frac{1} {2n}\biggr)-x
\biggr)^{-{n}/ {(n+1)}} \mathbh{1}_{[{1}/ 2 +{1}/ {(2(n+1))},
{1}/ 2 +{1} /{(2n)})}(x) \biggr],
\end{eqnarray}
where
%
%e13 #&#
%e14 #&#
\begin{eqnarray}
C_n&=&\int_{{1}/ 2 -{1}/ {(2n)}}^{{1}/ 2 -{1}/ {(2(n+1))}}
 \biggl(x-\biggl(
\frac{1} 2 -\frac{1} {2n}\biggr)\biggr)^{-{n}/ {(n+1)}} \,\mathrm{d}x
\nonumber
\\[-8pt]
\\[-8pt]
\nonumber
&&{}+ \int
_{{1}/ 2 +{1} /{(2(n+1))}}^{ {1}/ 2 +{1}/{(2n)}} \biggl(\biggl(\frac{1} 2 +
\frac{1} {2n}\biggr)-x\biggr)^{-{n}/ {(n+1)}} \,\mathrm{d}x,
\\
C&=&\Biggl(\sum_{n=1}^\infty\frac{1} {n^3}
\Biggr)^{-1}.
\end{eqnarray}

This density is quite similar to the one introduced in Counterexample
\ref{divergenza} and, in the same way we can prove that
$p\notin{\cal E}(1)$.

In order to see wether $p_t$ belongs to ${\cal E}(1)$ or not, we remark
that the same convergence problem arises whenever we integrate the
function $p^{1+\varepsilon}$ over any interval containing $\frac{1} 2$.
As a consequence, using the explicit formula \eqref{pt} for the
restriction $p_t$, we can prove that $p_t\notin{\cal E}(1)$, $\forall
 t\geq\frac{1} 2$. On the other hand,
for any $ t<\frac{1} 2$, we can find some $\varepsilon>0$ (depending on
$t$), such that $\E_{\mu}(p_t^{1+\varepsilon})<+\infty$.
Moreover, the condition $\E_{\mu}(p_t^{-\varepsilon})<+\infty$
is trivially satisfied, so that $p_t\in{\cal E}(1)$.

Finally, let us fix $t_0 <\frac{1} 2$ and define the density
\[
q(x)=p(x) \mathbh{1}_{[0,t_0]}(x)+\frac{1-F(t_0)}{1-t_0^\beta} \beta
 x^{\beta-1} \mathbh{1}_{(t_0,1]}(x),
\]
with $\beta>0$.
This function differs from the restriction $p_{t_0}$ only over
$(t_0,1]$, where $p_{t_0}$ is constant, while $q$ is proportional to a
beta density.

Using the arguments above and those in Example~\ref{beta} on beta
densities, we conclude that $q\in{\cal E}(1)$ (and $q_t\in{\cal
E}(1)$, $\forall 0\leq t<T$).
On the other hand, by construction, $p_t=\E_\mu(p|\mathcal F_t)=\E_\mu
(q|\mathcal F_t)$, $\forall t\leq t_0$.
\end{counterexample}

%re4.20 #&#
\begin{remark}
In a similar way, the density defined in Counterexample \ref
{divergenza} provides also a counterexample of Corollary~\ref{Exp1}.
\end{remark}

%s4.3 #&#
\subsection{Convexity and \texorpdfstring{$L^1(\mu)$}{L1(mu)}-closure}\label{geo}

Given $p\in{\cal P}$, we denote by ${\cal M}(p)$ the set of all
densities $q \in{\cal P}$ which are connected to $p$ by an open
mixture arc.
Moreover, let us recall that, by Theorem~\ref{main}, ${\cal E}(p)$
coincides with the set of all densities $q \in{\cal P}$ which are
connected to $p$ by an open exponential arc.

%pr4.21 #&#
\begin{proposition}\label{Prconv}
Let $p\in\mathcal P$. Then $\mathcal E(p)$ and $\mathcal M(p)$ are convex.
\end{proposition}

\begin{pf}
Note that for any $q, r\in\mathcal E(p)$, since $\mathcal
E(q)=\mathcal E(r)=\mathcal E(p)$, it is not restrictive to consider $r=p$.
Suppose $q\in\mathcal E(p)$, and consider $p(\lambda)=\lambda
p+(1-\lambda) q$ for any $\lambda\in[0,1]$.
We show that $p(\lambda)\in\mathcal E(p)$ by proving that   $\E_\mu
(p(\lambda)^\theta p^{1-\theta})<+\infty$  for   $\theta
\in(-\varepsilon,1+\varepsilon)$.

If $\theta\in(0,1)$, it follows by Jensen inequality
\[
\E_\mu\bigl(p(\lambda)^\theta p^{1-\theta}\bigr)=
\E_p\biggl(\biggl(\frac
{p(\lambda)}{p}\biggr)^\theta\biggr)
\le1.
\]
If $\theta\in(-\varepsilon,0)\cup(1,1+\varepsilon)$, by the convexity of
$x^\theta$ we have
\[
\E_\mu\bigl(p(\lambda)^\theta p^{1-\theta}\bigr)=
\E_p\biggl(\biggl(\frac
{\lambda p+(1-\lambda)q}{p} \biggr)^\theta\biggr)\le
\lambda+(1-\lambda )\E_p\biggl(\biggl(\frac{q}{p}
\biggr)^\theta\biggr)<+\infty,
\]
where the last inequality is due to $q\in\mathcal E(p)$.

On the other hand, one can easily see that $\mathcal M(p)$ is convex,
since the relation between open mixture arcs is an equivalence relation
and, thus, $q\in\mathcal M(p)$ implies $\mathcal M(p)=\mathcal M(q)$.
\end{pf}

In the following theorem we prove that the open mixture model $\mathcal
{M}(p)$ is $L^1(\mu)$-dense in the set of all densities $\mathcal P_{\ge
}$.

Since from Proposition~\ref{cor} $\mathcal{M}(p)\subseteq\mathcal
{E}(p)$, we deduce that $\mathcal{E}(p)$ is $L^1(\mu)$-dense in
$\mathcal P_{\ge}$. This result was proved in Imparato and Trivellato \cite{IT}, Theorem~19.1.

In the proof we use Scheff\'{e}'s theorem and, in particular, that if
$q,\{q_n\}_{n\ge1}$ are in $\mathcal P_{\ge}$ and such that $q_n\to
q$, as $n\to+\infty$, $\mu$-a.e., then
$q_n\to q$ in $L^1(\mu)$.

We remark that the $L^1(\mu)$ convergence when restricted to the set of
nonnegative densities $\mathcal P_{\ge}$, is equivalent to the
convergence in $\mu$-probability.

%th4.22 #&#
\begin{theorem}\label{Th1}
For any $p\in\mathcal P$ the open mixture model $\mathcal{M}(p)$ is
$L^1(\mu)$-dense in the nonnegative densities $\mathcal P_{\ge}$, that is
$\overline{\mathcal{M}(p)}=\mathcal P_{\ge}$, where the overline
denotes the closure in the $L^1(\mu)$-topology.
\end{theorem}
\begin{pf}
We show $\overline{\mathcal{M}(p)}=\mathcal P_{\ge}$,\vspace*{1pt} by checking the
double inclusions $\overline{\mathcal{M}(p)}\subseteq\mathcal P_{\ge}$
and $\overline{\mathcal{M}(p)}\supseteq\mathcal P_{\ge}$.

The first is straightforward, since by definition $\mathcal
M(p)\subseteq\mathcal P$ and, therefore, $\overline{\mathcal
{M}(p)}\subseteq\overline{\mathcal P}=\mathcal P_{\ge}$.
With regard to the second implication, we will start by showing that
any simple density $q\in\mathcal P_{\ge}$ is the $\mu$-a.e. limit of a
sequence of densities $q_n \in\mathcal{M}(p)$.

Let us fix a simple density $q\in\mathcal P_{\ge}$ and denote $\Omega
'=\operatorname{Supp}  q=\{\omega\dvt q(\omega)>0\}$ and $\Omega''=(\Omega')^c=\{
\omega\dvt q(\omega)=0\}$. Moreover, define the sequence $p_n$ by
\[
p_n=\frac{1}n \mathbh{1}_{(p<{1}/n)}+p \mathbh{1}_{({1}/n\le p\le n)}+n
\mathbh{1}_{(p>n)},
\]
and note that $p_n\to p$ pointwise.
For $n\ge1$, construct the sequence of densities in $\mathcal P$
\[
q_n=\cases{ %
\displaystyle\frac{qp}{p_n}
\frac{1}{c_n}, &\quad $\mbox{if } \omega\in\Omega',$
\vspace*{2pt}\cr
\displaystyle\frac{a_np}{c_n}, &\quad $\mbox{if } \omega\in\Omega'',$}
\]
where $a_n$ denotes a positive numerical sequence converging to 0 as
$n\to+\infty$ and $c_n=\int_{\Omega'} \frac{qp}{p_n}\,\mathrm{d}\mu+a_n\mathbb
{P}(\Omega'')$ the normalization constant.
As $q$ is a simple density whose support is $\Omega'$, there exist two
positive constants $a, A$ such that $a\le q(\omega)\le A$ if $\omega\in
\Omega'$.

Note that $q_n\in\mathcal M(p)$ by Theorem~\ref{main2}, since $\frac
{q_n}{p}$ is bounded from below and above, respectively, by two
positive constants:
\[
\frac{\min({a}/{n},a_n)}{c_n}\le\frac{q_n}{p}\le\frac{\max(An,a_n)}{c_n}.
\]
In order to check that $q_n\to q$, as $n\to+\infty$, $\mu$-a.e., it is
sufficient to prove that $c_n \to1$. In fact, since $\frac{p}{p_n}\le
\max(1,p)$ and the function
$q\max(1,p)\mathbh{1}_{\Omega'}\in L^1(\mu)$, we can apply the dominated
convergence theorem and find
\[
\lim_{n\to+\infty} c_n=\int_{\Omega'}
\lim_{n\to+\infty
} \frac{qp}{p_n}\,\mathrm{d}\mu=\int_{\Omega'}q\,\mathrm{d}
\mu=1.
\]
Thus, by Scheff\'{e}'s theorem the sequence $q_n\to q$ in $L^1(\mu)$.

Since any density $\mathcal P_{\ge}$ can be written as the limit of
simple densities in $L^1(\mu)$, we have proved that $\overline{\mathcal
{M}(p)}\supseteq\mathcal P_{\ge}$.
This concludes the proof of the theorem.
\end{pf}

The next corollary shows that the positive densities with finite
Kullback--Leibler divergence with respect to any $p\in\mathcal P$ is
$L^1(\mu)$-dense in the set of all densities $\mathcal P_{\ge}$. This
corresponds to the choice $\varphi(x)=x(\log(x))^+$ in the following result.

%co4.23 #&#
\begin{corollary}
Assume $\varphi\dvtx (0,+\infty)\to(0,+\infty)$ is a continuous function.
Then the set
\[
\mathcal P_\varphi=\biggl\{q\in\mathcal P\dvt \E_p\biggl(
\varphi\biggl(\frac
{q}{p}\biggr)\biggr)<+\infty\biggr\}
\]
is $L^1(\mu)$-dense in $\mathcal P_{\ge}$.
\end{corollary}
\begin{pf}
Let $q\in\mathcal P_\ge$. The result immediately follows from Theorem~\ref{Th1}, since the sequence $q_n\in\mathcal M(p)$ converging to $q$
in $L^1(\mu)$ is in $\mathcal P_\varphi$, that is
it satisfies $\E_p(\varphi(\frac{q_n}{p}))<+\infty$.
\end{pf}

%\appendixtitleon
%
%sA #&#
\begin{appendix}\label{app}
%%sB #&#
\section*{Appendix}

In this appendix we refer to Counterexample \ref{divergenza} and we
prove that both $D(q\|p)<+\infty$ and $D(p\|q)<+\infty$.

In fact,
\begin{eqnarray*}
D(q\|p)&=&\E_{\mu}(q\log q)=\int_0^1
q(x)\log q(x) \,\mathrm{d}x
\\
&=&\sum_{n=1}^\infty \int_{1 -{1} /{n}}^{1-{1}/ {(n+1)}}
\frac{C} {n^3 C_n} \biggl(x-\biggl(1 -\frac{1} {n}\biggr)
\biggr)^{-{n}/ {(n+1)}} \\
&&\hspace*{62pt}{}\times\log\biggl( \frac{C} {n^3 C_n} \biggl(x-\biggl(1 -
\frac{1} {n}\biggr) \biggr)^{-{n}/ {(n+1)}} \biggr) \,\mathrm{d}x
\\
&=&\sum_{n=1}^\infty \biggl[
\frac{C} {n^3 C_n}\log\biggl( \frac{C} {n^3 C_n}\biggr) \int_{1 -{1}/{n}}^{1-{1}/ {(n+1)}}
\biggl(x-\biggl(1 -\frac{1} {n}\biggr) \biggr)^{-{n}/ {(n+1)}} \,\mathrm{d}x
\\
&&\hspace*{18pt}{} +\frac{C} {n^3 C_n} \int_{1 -{1} /{n}}^{1-{1}/ {(n+1)}} \biggl(x-
\biggl(1 -\frac{1} {n}\biggr)\biggr)^{-{n} /{(n+1)}} \\
&&\hspace*{98pt}{}\times \log \biggl(\biggl(x-
\biggl(1 -\frac{1} {n}\biggr)\biggr)^{-{n}/ {(n+1)}} \biggr) \,\mathrm{d}x \biggr]
\\
&=&\sum_{n=1}^\infty \biggl[
\frac{C} {n^3}\log\biggl( \frac{C} {n^3 C_n}\biggr) -\frac{C} {n^3 C_n} \int
_{0}^{{1}/ {(n(n+1))}} \frac{n}{n+1} y^{-{n} /{(n+1)}} \log
y \,\mathrm{d}y \biggr]
\\
&=&\sum_{n=1}^\infty \biggl[
\frac{C} {n^3}\log\biggl( \frac{C} {n^3 C_n}\biggr) -\frac{C} {n^3}
\frac
{n}{n+1} \bigl(-\log\bigl(n(n+1)\bigr)-(n+1)\bigr)\biggr]
\\
&=&\sum_{n=1}^\infty \biggl[
\frac{C\log C} {n^3}-\frac{3C\log n}{n^3}+\frac{C} {n^3(n+1)} \log {\bigl(n(n+1)
\bigr)}-\frac{C} {n^3}\log{(n+1)}
\\
&&\hspace*{18pt}{} +\frac{C} {n^2(n+1)} \log\bigl(n(n+1)\bigr)+ \frac{C} {n^2} \biggr].
\end{eqnarray*}
Since the general term of the last series defines an infinitesimal
sequence of the same order as $\frac{1} {n^2}$, when $n\rightarrow\infty
$, we deduce that $D(q\|p)<+\infty$.

In the same way,
\begin{eqnarray*}
D(p\|q)&=&\E_{\mu}\biggl(\log\frac{1} q\biggr)=\int
_0^1 \log\frac{1}{
q(x)} \,\mathrm{d}x
\\
&=&\sum_{n=1}^\infty \int_{1 -{1}/ {n}}^{1-{1}/ {(n+1)}}
\log\biggl( \frac{C} {n^3 C_n} \biggl(x-\biggl(1 -\frac{1} {n}\biggr)
\biggr)^{-{n}/ {(n+1)}} \biggr)^{-1} \,\mathrm{d}x
\\
&=&\sum_{n=1}^\infty \int_{1 -{1} /{n}}^{1-{1}/ {(n+1)}}
\log\biggl( \frac{n^3 C_n} C \biggl(x-\biggl(1 -\frac{1} {n}\biggr)
\biggr)^{{n}/ {(n+1)}} \biggr) \,\mathrm{d}x
\\
&=&\sum_{n=1}^\infty \biggl[ \log\biggl(
\frac{n^3 C_n} C\biggr)\frac{1} {n(n+1)}+ \int_{0}^{{1}/ {(n(n+1))}}
\frac{n}{n+1} \log y \,\mathrm{d}y \biggr]
\\
&=&\sum_{n=1}^\infty \biggl[ \log\biggl(
\frac{n^3 C_n} C\biggr)\frac{1} {n(n+1)}+ \frac{1}{(n+1)^2} \bigl(-\log
\bigl(n(n+1)\bigr)-1\bigr) \biggr]
\\
&=&\sum_{n=1}^\infty \biggl[
\frac{3\log n}{n(n+1)} -\frac{\log C}{n(n+1)}+\frac{\log
(n+1)}{n(n+1)}-\frac{\log(n(n+1))}{n(n+1)^2}
\\
&&\hspace*{18pt}{} - \frac{\log(n(n+1))}{(n+1)^2} -\frac{1}{(n+1)^2} \biggr].
\end{eqnarray*}
As before, since the general term of the series defines an
infinitesimal sequence of the same order as $\frac{\log n} {n^2}$,
when $n\rightarrow\infty$, we get that $D(p\|q)<+\infty$.
\end{appendix}

\section*{Acknowledgements}  The authors are grateful to
Vincenzo Recupero for his fruitful suggestions on Counterexamples~\ref{divergenza}
and \ref{co4.19}.
M. Santacroce gratefully
acknowledges the hospitality of the Department of Mathematics
of the University of Texas at
Austin.

% imsref loaded by akundreckaite, 2015-04-29 08:50:41

%\begin{appendix}
%\section{}
%\end{appendix}

%\section*{Acknowledgements}

%\begin{supplement}%[id=suppA]
%\slink[doi]{10.3150/00-BEJXXXXSUPP} %[doi,text={...}] - jei reikia
%suskaldyti doi
%\sfilename{BEJ000\_supp.pdf}
%\end{supplement}

%\begin{thebibliography}{00}
%\bibitem{r1}
%\bibitem{r1}
%\end{thebibliography}

\printhistory
\end{document}